\documentstyle[11pt]{amsart}

\title[Schr\"odinger-type equations with variable exponent]
{Stationary waves of Schr\"odinger-type \\ equations with variable exponent}

\author{Du\v san D. Repov\v s}

\address{Faculty of Education, and
Faculty of Mathematics and Physics,
University of Ljubljana
\&
Institute of Mathematics, Physics and Mechanics, 1000 Ljubljana, Slovenia}
\email{dusan.repovs@@guest.arnes.si}

\newtheorem{theorem}{Theorem}[section]
\newtheorem{lemma}{Lemma}[section]
\newtheorem{proposition}{Proposition}[section]

\newtheorem{remark}{Remark}[section]

\def\RR{{\mathbb R}}

\def\phi{\varphi}

\def\ep{\varepsilon}

\def\phi{\varphi}

\def\di{\displaystyle}
\def\ri{\rightarrow}
\def\intom{\int_\Omega}

\def\ee{{\mathcal E}}

\def\w10{W^{1,p(x)}_0(\Omega)}
\begin{document}


\begin{abstract}
We are concerned with a class of nonlinear Schr\"odinger-type equations with a reaction term and a differential operator that involves a variable exponent. By using related variational methods, we establish several existence results.\\
\end{abstract}
\maketitle

\section{Introduction}
The topic of function spaces with variable exponents has undergone an impressive development over the last decades. It seems that the oldest references in this field are the works by  Orlicz \cite{orlicz} and  Nakano \cite{nakano}. This impressive revival is essentially connected with relevant applications to nonlinear partial differential equations  and fluid dynamics.

The field we are concerned with in this work  is of central interest, since the Schr\"odinger equation plays in quantum mechanics the same role as the Newton laws of conservation of energy play in classical mechanics. Schr\"odinger gave the classical derivation of this basic equation, based upon the analogy between mechanics and optics, and using the developments due to Louis de Broglie. The linear Schr\"odinger equation provides a thorough description of a particle in a non-relativistic setting. The structure of the nonlinear Schr\"odinger equation is much more complicated.
The nonlinear Schr\"odinger equation describes central phenomena arising in nonlinear optics, Bose-Einstein condensates, Heisenberg ferromagnets and magnons, plasma physics (the Kurihara superfluid film equation or Langmuir waves), condensed matter theory, etc. We refer to Ablowitz, Prinari and Trubatch \cite{abl} and Sulem \cite{sulem} for a modern overview, including relevant applications.

In this paper, we study a Schr\"odinger-type equation in the framework of function spaces with variable exponent. Problems of this type have been intensively studied in the last few years due to major applications to non-Newtonian electrorheological fluids (Halsey \cite{halsey}, Ruzicka \cite{ruz}) or image restoration (Chen, Levine and Rao \cite{chen}).  Our main purpose is to extend the study of Laplace-type operators to more general classes of equations with variable exponent and nonhomogeneous differential operator. We are concerned with Schr\"odinger-type equations that involve the differential operator $\mbox{div}\, (A(x,|\nabla u|)\nabla u)$ and power-type nonlinearities with exponent variable. 
The abstract setting in the present paper corresponds to very general differential operators that include the usual $p(x)$-Laplace or $p(x)$-mean curvature operators, resp.
$$\mbox{div}\,(|\nabla u|^{p(x)-2}\nabla u),\quad \mbox{div}\,(p(x)|\nabla u|^{p(x)-2}\nabla u)\quad\mbox{and}\quad \mbox{div}\,\left((1+|\nabla u|^2)^{(p(x)-2)/2}\nabla u\right).$$
An excellent reference for the most significant mathematical
methods employed in this paper is the 
book by Ciarlet \cite{ciarlet}.

\section{A brief review on variable exponent Lebesgue-Sobolev spaces}
Throughout this paper we assume that $\Omega$ is a bounded open set in $\RR^N$ with smooth boundary.

In this section, we recall some definitions and basic properties of the variable exponent Lebesgue space $L^{p(x)}(\Omega)$ and $W_0^{1,p(x)}(\Omega)$. Roughly speaking,  Lebesgue and Sobolev spaces with variable exponent are functional spaces of Lebesgue's and Sobolev's type in which different space directions have different roles. The variable exponent Lebesgue space $L^{p(x)}(\Omega)$ is a special case of Orlicz-Musielak spaces treated by Musielak \cite{musielak}.

 Denote
$$
C_+(\overline{\Omega})=\{h;~h\in C(\overline{\Omega}),~h(x)>1,~\mbox{for all}\ x\in \overline{\Omega}\}.
$$
For any $h\in C_+(\overline{\Omega})$, we define
$$
h^+=\max\{h(x);~x \in\overline{\Omega}\},\hspace{0.3cm}h^-=\min\{h(x);~x \in\overline{\Omega}\}.
$$
For any $p\in C_+(\overline{\Omega})$, we define the $\emph{variable exponent Lebesgue space}$
$$
L^{p(x)}(\Omega)=\Big\{u:\Omega\rightarrow\RR;~ u\textrm{ is  measurable and } \int_{\Omega}|u(x)|^{p(x)}dx<\infty\Big\},
$$
endowed with the  $\emph{Luxemburg norm}$
$$
|u|_{L^{p(x)}(\Omega)}=|u|_{p(x )}=\inf \Big\{\mu>0;~\int_{\Omega}\left|\frac{u(x)}{\mu}\right|^{p(x)}dx\leq 1\Big\}.
$$
Then $(L^{p(x )}(\Omega),|.|_{p(x )})$ is a Banach space, cf. \cite{Kov}.

As established in \cite{Ed1},
 $(L^{p(x )}(\Omega),|.|_{p(x )})$ is a separable, uniformly convex Banach space and its dual space is $L^{q(x )}(\Omega)$, where $\frac{1}{p(x )}+\frac{1}{q(x )}=1.$ Moreover, for all  $u\in L^{p(x )}(\Omega)$ and $v\in L^{q(x )}(\Omega)$, we have the H\"older inequality
$$
\Big{|}\int_{\Omega}uvdx\Big{|}\leq \Big(\frac{1}{p^{-}}+
\frac{1}{q^{-}}\Big)|u|_{p(x )}|v|_{q(x )}\,.
$$
If $p_1(x ),~p_2(x )\in C_{+}(\overline{\Omega})$ and $p_1(x )\leq p_2(x )$ for all $x\in \overline{\Omega}$, then $L^{p_2(x )}(\Omega)\hookrightarrow L^{p_1(x )}(\Omega)$ and the embedding is continuous.

An important role in manipulating the  generalized Lebesgue space is played by the $p(x )$-modular of the $L^{p(x )}(\Omega)$ space, which is the mapping $\rho_{p(x )}:L^{p(x )}(\Omega)\rightarrow \RR$ defined by
$$
\rho_{p(x )}(u)=\int_{\Omega}|u|^{p(x)}dx.
$$
\begin{proposition}[See \cite{Fan4}]\label{prop2.2}
For $u\in L^{p(x )}(\Omega)$ and $u_n\subset L^{p(x )}(\Omega)$, we have\\
(1)\ $|u|_{p(x )}<1$ \ (respectively $=1; >1$)$\Longleftrightarrow \rho_{p(x )}(u)<1$ \ (respectively $=1; >1$);\\
(2)\ for $u\neq 0$, $|u|_{p(x )}=\lambda \Longleftrightarrow \rho_{p(x )}(\frac{u}{\lambda})=1$;\\
(3)\ if $|u|_{p(x )}>1$, then $|u|_{p(x )}^{p^{-}}\leq \rho_{p(x )}(u)\leq |u|_{p(x )}^{p^{+}}$;\\
(4)\ if $|u|_{p(x )}<1$, then $|u|_{p(x )}^{p^{+}}\leq \rho_{p(x )}(u)\leq|u|_{p(x )}^{p^{-}}$;\\
(5)\ $|u_n-u|_{p(x )}\rightarrow0$ (respectively $\rightarrow \infty$)$\Longleftrightarrow \rho_{p(x )}(u_n-u)\rightarrow 0$ (respectively $\rightarrow \infty$),
since $p^+<\infty$.
\end{proposition}

The  Sobolev space with variable exponent $W^{1,p(x)}(\Omega)$ is defined by
$$
W^{1,p(x)}(\Omega)=\{u\in L^{p(x)}(\Omega);~\partial_{x_i}u\in L^{p(x)}(\Omega),~i\in \{1, \dots ,N\}\}.
$$
If equipped with the norm
$$
\|u\|_{W^{1,p(x)}(\Omega)}=|u|_{L^{p(x)}(\Omega)}+\sum_{i=1}^N|\partial_{x_i}u|_{L^{p(x)}(\Omega)},
$$
then ($W^{1,p(x)}(\Omega),\|\cdot\|_{W^{1,p(x)}(\Omega)}$) is a separable and reflexive Banach space (see \cite[Theorem 1.3]{Kov}).

As observed by Zhikov \cite{zh1} in relationship with the {\it Lavrentiev phenomenon}, it is possible that minimizers of certain variational integrals  are not smooth. The following {\it log-H\"older condition} was first used in the variable exponent context by Zhikov \cite{zh2}. More precisely, we say that a function $h:\Omega\ri\RR$ is log-H\"older continuous on $\Omega$ if there exists $C>0$ such that
\begin{equation}\label{loghol}|h(x)-h(y)|\leq-\frac{C}{\log |x-y|}\qquad\mbox{for all $x,y\in\Omega$, $|x-y|\leq\frac 12$.}\end{equation}
As established in \cite{Ed1} (see also \cite[Theorem 9.1.8]{DHHR}), if $\Omega$ is bounded with Lipschitz boundary and $p$ is log-H\"older continuous, then $C^\infty(\overline\Omega)$ is dense in $W^{1,p(x)}(\Omega)$.
We point out  that though log-H\"older continuity of $p(x)$ is sufficient to imply the density of test functions in $W^{1,p(x)}$, this condition is far from being necessary. For instance, Edmunds and R\'akosn\'{\i}k \cite{Ed0} derived the same conclusion under a local monotony condition on $p$.

Let $p$ be log-H\"older continuous. The Sobolev space  $W_0^{1,p(x )}(\Omega)$  with zero boundary values is the closure of the set of $W^{1,p(x)}(\Omega)$-functions with compact support. Furthermore, if $p$ is bounded then $W_0^{1,p(x )}(\Omega)$ is the closure of $C^\infty_0(\Omega)$
in the space $W^{1,p(x)}(\Omega)$, see \cite[Proposition 11.2.3]{DHHR}.
The norm $||u||=\sum_{i=1}^N|\partial_{x_i}u|_{p(x)}$ is an equivalent norm in $W_0^{1,p(x )}(\Omega)$ (see \cite{MPR}). Hence $W_0^{1,p(x )}(\Omega)$ is a separable and reflexive Banach space. Note that when $s\in C_+(\overline{\Omega})$ and $s(x)<p^*(x)$ for all $x\in \overline{\Omega}$, where $p^*(x)=\frac{Np(x)}{N-p(x)}$ if $p(x)<N$ and $p^*(x)=\infty$ if $p(x)\geq N$, then the embedding $W_0^{1,p(x )}(\Omega)\hookrightarrow L^{s(x )}(\Omega)$ is compact.

\section{Main results}
In this paper we study the following nonlinear problem
\begin{equation}\label{problem}
\left\{\begin{array}{lll}
&\di -\mbox{div}\, (A(x,|\nabla u|)\nabla u)= f(x,u)&\quad\mbox{if $x\in\Omega$}\\
&\di u=0&\quad\mbox{if $x\in\partial\Omega$}.\end{array}\right.
\end{equation}

Problems of this type are motivated by models in mathematical physics (see Reed and Simon \cite{reed} and Strauss \cite{strauss}), where certain stationary waves in nonlinear Schr\"odinger or Klein-Gordon equations can be reduced to this form. Equation \eqref{problem} arises in the study of the Schr\"odinger-type equation
$$iv_t-\mbox{div}\, (A(x,|\nabla v|)\nabla v)= f(x,v)$$
when looking for {\it standing waves}, that is, solutions of the type $v(x,t)=e^{-ict}u(x)$, where $c$ is a real constant. This problem has a central role in quantum mechanics for the study of particles of stochastic fields modelled by L\'evy processes. A path integral over the L\'evy flights paths and a nonlinear Schr\"odinger equation is formulated by Laskin \cite{laskin} from the idea of Feynman and Hibbs's paths integrals.

\smallskip
Fix $p\in C_+(\overline{\Omega})$ and suppose that $p$ is log-H\"older continuous.

We assume that $A:\Omega\times [0,\infty)\ri [0,\infty)$ satisfies the following hypotheses:

\smallskip\noindent
(A1) the mapping
$\Omega\ni x\longmapsto A(x,s)$ is measurable for all $s\geq 0$ and the mapping
$(0,\infty)\ni s\longmapsto A(x,s)$ is absolutely continuous for a.a. $x\in \Omega$;

\smallskip\noindent
(A2) there exist $a_1\in L^{p'(x)}(\Omega)$ and $a_2>0$ such that
$$ |A(x,|z|)z|\leq a_1(x)+a_2\, |z|^{p(x)-1}\quad \mbox{a.e. $x\in\Omega$, all $z\in\RR^N$;}$$
(A3) there exists $a_3>0$ such that for a.a. $x\in\Omega$ and for all $s>0$
$$\min\{A(x,s),A(x,s)+s\partial_s A(x,s)\}\geq a_3\,\min\{1,s^{p(x)-2}\}\,.$$
(A4) we have $\di t^2A(x,t)\leq p_+\int_0^{t}sA(x,s)ds$ for a.a. $x\in\Omega$ and for all $t\in (0,\infty)$.

\smallskip
We assume that the nonlinear term $f:\Omega\times\RR\ri\RR$ is a Carath\'eodory function satisfying the following conditions:

\smallskip
\noindent (f1) there exists $C>0$ such that $|f(x,t)|\leq C\,|t|^{q(x)-1}$ for a.a. $x\in\Omega$ and for all $t\in\RR$, where $q\in C_+(\overline\Omega)$ and $\max_{x\in\overline\Omega}p(x)<\min_{x\in\overline\Omega}q(x)$, $q(x)<p^*(x)$ for all $x\in\overline\Omega$;

\smallskip
\noindent (f2) there exist $\mu>p_+$ and $R>0$ such that
$$0<\mu F(x,t)\leq tf(x,t)\qquad\mbox{for all $x\in\Omega$ and for all $t\geq R$},$$
where $F(x,t):=\int_0^tf(x,s)ds$;

\smallskip
\noindent
(f3) $\lim_{t\ri 0}f(x,t)/|t|^{p_+-1}=0$ uniformly for $x\in\Omega$.

\smallskip
We say that $u$ is a solution of problem \eqref{problem} if $u\in W_0^{1,p(x)}(\Omega)\setminus\{0\}$ and
$$\int_\Omega A(x,|\nabla u|)\nabla u\cdot\nabla v dx=\lambda\int_\Omega f(x,u)vdx\quad\mbox{for all}\ v\in W_0^{1,p(x)}(\Omega).$$

The first result in this paper establishes the existence of solutions to problem \eqref{problem} under the above hypotheses. This result extends previous existence properties
to a very large class of nonlinear differential operators. The proof combines variational arguments and related energy estimates. A key role is played by the Mountain Pass Theorem, see Ambrosetti and Rabinowitz \cite{amra}.

\begin{theorem}\label{t1}
Assume that $A$ and $f$ satisfy conditions (A1)--(A4) and (f1)--(f3). Then problem \eqref{problem} has at least one solution.
\end{theorem}

Next, we are concerned with the existence of multiple high-energy solutions of problem \eqref{problem}. For this purpose, we consider the associated energy functional $\ee :\w10\ri\RR$ defined by
$$\ee (u)=\intom\left(\int_0^{|\nabla u(x)|}sA(x,s)ds\right)dx-\intom F(x,u)dx.$$

Assuming that $f(x,\cdot )$ is odd, we prove that problem \eqref{problem} has a sequence of solutions with higher and higher energies. The statement of this result is the following.

\begin{theorem}\label{t2}
Suppose that hypotheses (A1)--(A4) and (f1)-(f2) are fulfilled and $f(x,-t)=-f(x,t)$ for a.a. $x\in\Omega$ and all $t\in\RR$. Then problem \eqref{problem} admits a sequence of solutions $(u_n)$ such that $\ee (u_n)\ri+\infty$ as $n\ri\infty$.
\end{theorem}

A central
role  in the proof of Theorem \ref{t2} is played by the Fountain Theorem, which is due to
Bartsch \cite{ref91}. This result is nicely presented in Willem \cite{Will} by using the quantitative
deformation lemma. We also point out that the dual version of the Fountain Theorem is due to Bartsch and Willem, see \cite{Will}. It should be noted that the Palais-Smale condition
plays an important role for these theorems and their applications.

A related question concerns the minimization problem corresponding to the associated Rayleigh quotient, namely
$$\lambda_1:=\inf_{u\in\w10\setminus\{0\}}\frac{\di\intom\left(\int_0^{|\nabla u(x)}sA(x,s)ds\right)dx}{\di\intom F(x,u)dx}\,.$$

A natural assumption, in order to avoid degeneracies, is the following.

\smallskip
\noindent (f4) for all $(x,u)\in\Omega\times\RR$ we have $uf(x,u)\geq 0$ and $f\not\equiv 0$.

\smallskip
In the case corresponding to the Laplace or $p$-Laplace operators we have $\lambda_1>0$. However, it may happen that $\lambda_1=0$, see \cite{radpams}. The next result establishes a sufficient condition such that $\lambda_1$ is positive. This condition takes into account the growth of the variable potentials $p(\cdot)$ and $q(\dot)$ and is the following.

\smallskip
\noindent (A5) we have $2(q_+-q_-)<p_-$.

\begin{theorem}\label{t3}
Assume that hypotheses (A1)--(A5), (f1) and (f4) are fulfilled. Then 
\begin{equation}\label{minipro}
\lambda_1:=\inf_{u\in\w10\setminus\{0\}}\frac{\di\intom\left(\int_0^{|\nabla u(x)}sA(x,s)ds\right)dx}{\di\intom F(x,u)dx}>0.
\end{equation}
Moreover, for all $\lambda\geq\lambda_1$ the nonlinear problem
\begin{equation}\label{problem1}
\left\{\begin{array}{lll}
&\di -{\mathrm div}\, (A(x,|\nabla u|)\nabla u)= \lambda f(x,u)&\quad\mbox{if $x\in\Omega$}\\
&\di u=0&\quad\mbox{if $x\in\partial\Omega$}\end{array}\right.
\end{equation}
has a nontrivial solution.
\end{theorem}

In particular, Theorem \ref{t3} establishes a concentration property at infinity for the spectrum
problem $Su=\lambda Tu$, where $Su:=-\mbox{div}\, (A(x,|\nabla u|)\nabla u)$ and $Tu:=f(x,u)$.

\section{Auxiliary properties}
From now on we assume that $A$ and $f$ satisfy the hypotheses in the previous section.

We first establish that the energy $\ee$ is of class $C^1$ on $\w10$ and for all $u,v\in\w10$
\begin{equation}\label{deriv}
\ee '(u)(v)=\intom A(x,|\nabla u|)\nabla u\nabla vdx-\intom f(x,u)vdx.
\end{equation}
Relation \eqref{deriv} follows if we prove that
\begin{equation}\label{deriv1}
\ee'_0(u)(v)=\intom A(x,|\nabla u|)\nabla u\cdot\nabla vdx,
\end{equation}
where
$$\ee_0 (u)=\intom\left(\int_0^{|\nabla u(x)|}sA(x,s)ds\right)dx.$$

Fix $u,v\in\w10$ and  $\lambda\in\RR$ close to zero, say
$0<|\lambda|<1$. Thus, by hypothesis (A2) combined with the mean value theorem, there exists
$\theta$ between 0 and $\lambda$ such that
\begin{equation}\label{l1}\begin{array}{ll}
\di&\di\frac{\di\int_0^{|\nabla u(x)+\lambda \nabla
v(x)|}A(x,s)sds-\int_0^{|\nabla u(x)|}A(x,s)sds}{\lambda}=\\ & \ \\
&\di A(x,|\nabla u(x)+\theta\nabla v(x)|)|\nabla u(x)+\theta\nabla
v(x)|\, |\nabla v(x)|\leq\\ & \ \\
&\di (a_1(x)+a_2|\nabla u(x)+\theta\nabla v(x)|^{p(x)-1})\, |\nabla
v(x)|\leq\\ & \ \\
&\di \leq \left(a_1(x)+a_2(|\nabla u(x)|+|\nabla v(x)|)^{p(x)-1}\right|)\,|\nabla
v(x)|.\end{array}\end{equation}
Next, we show that the right-hand side of relation \eqref{l1} is in
$L^1(\Omega)$. Indeed, by  H\"older's inequality, we have
$$\begin{array}{ll}
\di&\di\intom \left(a_1(x)+a_2(|\nabla u(x)|+|\nabla v(x)|)^{p(x)-1}\right)\,|\nabla
v(x)|dx\leq \\&\di 2|a_1|_{p'(\cdot)}\,|\nabla
v|_{p(\cdot)}+2a_2|(|\nabla u|+|\nabla
v|)^{p(x)-1}|_{p'(\cdot)}|\nabla v|_{p(\cdot)}\leq \\
&\di 2|a_1|_{p'(\cdot)}\,|\nabla v|_{p(\cdot)}+
C\left(\intom(|\nabla u(x)|^{p(x)}+|\nabla v(x)|^{p(x)})dx
\right)^{1/p'_-}\,|\nabla v |_{p(\cdot)}\,.\end{array}$$

Taking $\lambda\ri 0$, we also have $\theta\ri 0$. Thus, by the
Lebesgue dominated convergence theorem,
$$\begin{array}{ll}
\di \lim_{\lambda\ri 0}\frac{\ee_0 (u+\lambda v)-\ee_0
(u)}{\lambda}&\di =\intom\lim_{\theta\ri 0}A(x,|\nabla u+\theta\nabla
v|)(\nabla u+\theta\nabla v)\cdot\nabla vdx\\
&\di =\intom A(x,|\nabla u|)\nabla u\cdot\nabla vdx,\end{array}$$
which proves our claim \eqref{deriv1}.

Next, we consider the following operators:

(i) the gradient operator $\nabla :\w10 \ri L^{p(x)}(\Omega,\RR^N)$;

(ii) the Nemytskii operator $T:L^{p(x)}(\Omega,\RR^N)\ri
L^{p'(x)}(\Omega,\RR^N)$ defined by
$Tu(x)=A(x,| u(x)|)u(x)$ for all $u\in L^{p(x)}(\Omega,\RR^N)$;

(iii) the linear operator $L:L^{p'(x)}(\Omega,\RR^N)\ri
W^{-1,p'(x)}(\Omega)$ defined by $$Lu(v)=\intom u(x)\cdot\nabla v(x)dx\quad
\mbox{for all $u\in L^{p'(x)}(\Omega,\RR^N)$ and
$v\in \w10 $}.$$

These operators are continuous and $\ee_0'=L\circ T\circ\nabla$,
hence $\ee_0$ is of class $C^1$.

\smallskip
The next result is a counter-part of formula (2.2) in Simon \cite{simon}. More precisely, in \cite{simon} it is established that for all $\xi$, $\zeta\in\RR^N$
$$
|\xi-\zeta|^{p}\leq\begin{cases} c(|\xi|^{p-2}\xi- |\zeta|^{p-2}\zeta
)(\xi-\zeta)\qquad&\mbox{for }\phantom{1<\,}p\geq 2;\\
c\langle|\xi|^{p-2}\xi-|\eta|^{p-2}\eta,\xi-\eta \rangle^{p/2}
\left(|\xi|^p+|\eta|^p\right)^{(2-p)/2}\qquad&\mbox{for }1<p<2,\end{cases}$$
where $c$ is a positive constant.

As above, 
 we distinguish between the {\it singular} case corresponding to $1<p(x)<2$ and the {\it degenerate} case, which corresponds to $p(x)>2$. The version of the above inequalities for variable exponents is the following.

\begin{lemma}\label{lem1}
Assume that hypotheses (A1) and (A3) are fulfilled. Then there exists a positive constant $C$ such that for all $\xi,\zeta\in\RR^N$ with $(\xi,\zeta)\not=(0,0)$, we have:

(i) $\left(A(x,|\xi|)\xi-A(x,|\zeta|)\zeta\right)\cdot(\xi-\zeta)\geq C_1\,|\xi-\zeta|^{p(x)}$ for all $x\in\Omega$ with $p(x)\geq 2$;

(ii) $\left(A(x,|\xi|)\xi-A(x,|\zeta|)\zeta\right)\cdot (\xi-\zeta) \geq C_2\,|\xi-\zeta|^{2}\min\{1, (|\xi|+|\zeta|)^{p(x)-2}\}$ for all $x\in\Omega$ with $1<p(x)< 2$.
\end{lemma}

{\it Proof.} We have
$$\begin{array}{ll}
\di \left(A(x,|\xi|)\xi-A(x,|\zeta|)\zeta\right)\cdot (\xi-\zeta) &=\di \sum_{i=1}^N
\left(A(x,|\xi|)\xi_i-A(x,|\zeta|)\zeta_i\right) (\xi_i-\zeta_i) \\
& =\di \sum_{i=1}^N\left(\varphi_i(x,\xi)-\varphi(x,\zeta)\right) (\xi_i-\zeta_i),\end{array}$$
where
$\varphi_i(x,w):= A(x,|w|)w_i$ for all $w\in\RR^N$. But
$$\varphi_i(x,\xi)-\varphi(x,\zeta)=\sum_{j=1}^N\int_0^1\frac{\partial\varphi_i(x,z)}{\partial z_j}(\xi_j-\zeta_j)dt,$$
where $z=\zeta+t(\xi-\zeta)$. Therefore
\begin{equation}\label{sstar}
\left(A(x,|\xi|)\xi-A(x,|\zeta|)\zeta\right)\cdot (\xi-\zeta) =
\sum_{i,j=1}^N\int_0^1\frac{\partial\varphi_i(x,z)}{\partial z_j}(\xi_i-\zeta_i)(\xi_j-\zeta_j)dt.\end{equation}
Fix $z,w\in\RR^N\setminus\{0\}$. We observe that
\begin{equation}\label{sstar1}\begin{array}{ll}
\di \sum_{i,j=1}^N\frac{\partial\varphi_i(x,z)}{\partial z_j}w_iw_j&=\di A(x,|z|)|w|^2+
\frac{1}{|z|}\,A_s(x,|z|)\, (z\cdot w)^2\\
&\di = |w|^2\left(A(x,|z|)+|z|A_s(x,|z|)\, ( \frac{z}{|z|}\cdot \frac{w}{|w|})^2\right).\end{array}\end{equation}

(i) Assuming that $p(x)\geq 2$, hypothesis (A3) and relation \eqref{sstar1} yield
$$\sum_{i,j=1}^N\frac{\partial\varphi_i(x,\xi)}{\partial \xi_j}\zeta_i\zeta_j\geq a_3\,|\xi|^{p(x)-2}\, |\zeta|^2\,.$$
Returning to relation \eqref{sstar} we deduce that for all $t>0$
\begin{equation}\label{finale} \left(A(x,|\xi|)\xi-A(x,|\zeta|)\zeta\right)\cdot (\xi-\zeta)\geq a_3\int_0^1|\zeta +t(\xi-\zeta)|^{p(x)-2}\, |\xi-\zeta|^2\,dt.\end{equation}
Choosing $t=1/3$ we deduce that
$$|\zeta +t(\xi-\zeta)|\geq\max\{|\xi| ,|\zeta|\}-\frac{|\xi -\zeta|}{3}\geq \frac{|\xi -\zeta|}{3}\,.$$
Thus, by \eqref{finale}, we conclude that there exists $C_1>0$ such that
$$\left(A(x,|\xi|)\xi-A(x,|\zeta|)\zeta\right)\cdot (\xi-\zeta)\geq C_1\, |\xi -\zeta|^{p(x)}.$$

(ii) Assume that $1<p(x)<2$. Using assumption (A3) we have
$$\sum_{i,j=1}^N\frac{\partial\varphi_i(x,\xi)}{\partial \xi_j}\zeta_i\zeta_j\geq a_3\,\min\{1,|\xi|^{p(x)-2}\}\, |\zeta|^2\,.$$
It follows that
$$\begin{array}{ll}
&\di \left(A(x,|\xi|)\xi-A(x,|\zeta|)\zeta\right)\cdot (\xi-\zeta)\di=\\
& \di\sum_{i,j=1}^N\int_0^1\frac{\partial\varphi_i}{\partial z_j}(x,\zeta +t(\xi-\zeta))(\xi_i-\zeta_i)
(\xi_j-\zeta_j)dt \geq\\
&\di a_3\int_0^1\min\{1, |\zeta +t(\xi-\zeta)|^{p(x)-2}\}\, |\xi-\zeta|^2dt\geq\\
&\di C\, |\xi-\zeta|^2\min\{1, (|\xi| +|\zeta|)^{p(x)-2}\}.\end{array}$$
This completes the proof.
\qed

\begin{remark}\label{rem1}
Taking into account the expression of $\ee_0'$  in \eqref{deriv1}, Lemma \ref{lem1} implies that the operator $\ee_0': \w10\ri W^{-1,p'(x)}(\Omega)$ is strictly monotone. This fact is in accordance with the hypotheses imposed in \cite[Section 13.4]{DHHR}.
\end{remark}

\section{Proof of Theorem \ref{t1}}
By the results in the previous section, in order to find a solution of problem \eqref{problem}, it is enough to show that the energy functional $\ee$ has a nontrivial critical point. For this purpose we first check that $\ee$ has the Mountain Pass geometry and then we show that $\ee$ satisfies the Palais-Smale condition.

By hypotheses (f1) and (f3), for all $\ep>0$ there exists $C_\ep>0$ such that for all $(x,u)\in\Omega\times\RR$
$$|F(x,u)|\leq\ep\, |u|^{p_+}+C_\ep\, |u|^{q(x)}\,.$$
Thus, by (A3) and (A4)
$$\begin{array}{ll}
\di \ee (u) &\di \geq C\, |\nabla u|^{p_+}_{p(\cdot)}-\intom \left(\ep\, |u|^{p_+}+C_\ep\, |u|^{q(x)}\right)dx\\
&\di \geq C_1\, \| u|^{p_+}_{\w10}-\ep\, \| u|^{p_+}_{\w10}-C_2(\ep)\, |u|^{q_-}_{q(\cdot)}\,.\end{array}$$
By (f1) we know that $p_+<q_-$. This implies that there exists $\ep>0$ small enough and there are positive numbers $r$ and $\rho$ such that $\ee (u)\geq\rho$ for all $u\in\w10$ with $\| u|_{\w10}=r$.

\smallskip
We now check the second geometric assumption of the Mountain Pass Theorem, namely the existence of a ``valley". Fix $\varphi\in \w10\setminus\{0\}$ and $t>0$. We have
\begin{equation}\label{geo1}
\ee (t\varphi)=\intom\left(\int_0^{t|\nabla\varphi (x)|}sA(x,s)ds\right)dx-\intom F(x,t\varphi)dx.
\end{equation}
Hypothesis (f3) implies that there are $\alpha>p_+$  and positive constants $A$, $B$ such that for all $(x,u)\in\Omega\times\RR$
\begin{equation}\label{geo2}
F(x,u)\geq A\,|u|^\alpha -B.
\end{equation}
For fixed $x\in\Omega$ and $z\in\RR^N$, consider the differentiable function 
$$h(t):=\int_0^{t|z|}sA(x,s)ds\qquad t>0.$$
It follows that $h'(t)=t\,|z|^2A(x,t|z|)$. Thus, by (A4),
$$h'(t)=\frac 1t\, (t|z|)^2A(x,t|z|)\leq\frac{p_+}{t}\int_0^{t|z|}sA(x,s)ds=\frac{p_+}{t}\, h(t).$$
By integration we deduce that for all $t>0$
\begin{equation}\label{geo3}h(t)=\int_0^{t|z|}sA(x,s)ds\leq Ct^{p_+},\end{equation}
where
$$C=C(x,z)=\int_0^{|z|}sA(x,s)ds.$$
Using estimates \eqref{geo2} and \eqref{geo3}, relation \eqref{geo1} yields
\begin{equation}\label{geo4}
\ee (t\varphi)\leq C(\varphi)t^{p_+}-C_1t^\alpha +C_2,
\end{equation}
where $C(\varphi)$, $C_1$, $C_2$ are positive constants. Since $\alpha>p_+$, relation \eqref{geo4}
shows that $\ee (t\varphi)<0$ for $t$ large enough.

\smallskip
Next, we show that $\ee$ satisfies the Palais-Smale condition. For this purpose we need the following auxiliary result, which extends a classical property in functional analysis.

\begin{lemma}\label{weakco}
Let $(u_n)$ be a sequence in $\w10$ that converges weakly to $u$ and such that 
$$\limsup_{n\ri\infty} \left(\ee_0'(u_n)-\ee_0'(u)\right)(u_n-u)\leq 0.$$
Then $(u_n)$ converges strongly in $\w10$.
\end{lemma}

{\it Proof.} We already know (see Remark \ref{rem1}) that $\ee_0'$ is a monotone operator, hence
$$\left(\ee_0'(u_n)-\ee_0'(u)\right)(u_n-u)\geq 0.$$
Therefore
\begin{equation}\label{egal}
\lim_{n\ri\infty} \left(\ee_0'(u_n)-\ee_0'(u)\right)(u_n-u)= 0.
\end{equation}
But
$$\left(\ee_0'(u_n)-\ee_0'(u)\right)(u_n-u)= \intom\left( A(x,|\nabla u_n|)\nabla u_n-
A(x,|\nabla u|)\nabla u\right)\cdot(\nabla u_n-\nabla u)\,dx.$$

We claim that
\begin{equation}\label{cl1}
\intom |\nabla u_n-\nabla u|^{p(x)}dx\ri 0\quad\mbox{as $n\ri\infty$},
\end{equation}
which implies that $(u_n)$ converges strongly to $u$ in $\w10$.

We have
$$\intom |\nabla u_n-\nabla u|^{p(x)}dx=\int_{\Omega_+} |\nabla u_n-\nabla u|^{p(x)}dx+\int_{\Omega_-} |\nabla u_n-\nabla u|^{p(x)}dx,$$
where
$$\Omega_+:=\{x\in\Omega;\ p(x)\geq 2\}\quad \Omega_-:=\{x\in\Omega;\ 1<p(x)< 2\}.$$

Using Lemma \ref{lem1} we obtain
\begin{equation}\label{eu1}\begin{array}{ll}
&\di\int_{\Omega_+}\left( A(x,|\nabla u_n|)\nabla u_n-A(x,|\nabla u|)\nabla u\right)\cdot(\nabla u_n-\nabla u)dx\geq \\ &\di C_1\int_{\Omega_+}|\nabla u_n-\nabla u|^{p(x)}dx\end{array}\end{equation} and
\begin{equation}\label{eu2}\begin{array}{ll}
&\di\int_{\Omega_-}\left( A(x,|\nabla u_n|)\nabla u_n-A(x,|\nabla u|)\nabla u\right)\cdot(\nabla u_n-\nabla u)dx\geq\\ &\di C_2\int_{\Omega_-}\min\{1,(|\nabla u_n|+|\nabla u|)^{p(x)-2}\}\,|\nabla u_n-\nabla u|^{2}dx.\end{array}\end{equation}
Applying H\"older's inequality we obtain
\begin{equation}\label{eu3}
\int_{\Omega_-} |\nabla u_n-\nabla u|^{p(x)}dx\leq C\left(\int_{\Omega_-}\min\{1,(|\nabla u_n|+|\nabla u|)^{p(x)-2}\}\,|\nabla u_n-\nabla u|^{2}dx \right)^\gamma ,
\end{equation}
where $\gamma$ is a positive constant.

Relations \eqref{eu1}, \eqref{eu2} and \eqref{eu3} imply our claim \eqref{cl1}. This concludes the proof of Lemma. \qed

\smallskip
Returning to the proof of Theorem \ref{t1}, it remains to show that $\ee$ satisfies the Palais-Smale condition. For this purpose, let $(u_n)$ be a sequence in $\w10$ such that $\ee (u_n)\ri c\in\RR$ and $\ee' (u_n)\ri 0$ in $W^{-1,p'(x)}(\Omega)$. We first prove that \begin{equation}\label{roco}
\mbox{$(u_n)$ is bounded in $\w10$}.\end{equation} Using hypotheses (A4) and (f2) we have as $n\ri\infty$
$$\begin{array}{ll}
\di O(1)+o(\|u_n\|)&\di =\ee (u_n)-\frac 1\mu \ee '(u_n)(u_n)\\
&\di =\intom\left( \int_0^{|\nabla u_n|}sA(x,s)ds-\frac 1\mu\, A(x,|\nabla u_n|)|\nabla u_n|^2\right)dx\\ &\di +\intom (f(x,u_n)u_n-F(x,u_n))dx\\
&\di\geq \left(1-\frac{p_+}{\mu}\right)\intom\left( \int_0^{|\nabla u_n|}sA(x,s)ds\right)dx\\ &\di +
\intom (f(x,u_n)u_n-F(x,u_n))dx\\
&\di\geq \left(1-\frac{p_+}{\mu}\right)\intom\left( \int_0^{|\nabla u_n|}sA(x,s)ds\right)dx+O(1).\end{array}$$
Since $\mu>p_+$, we obtain
$$\intom\left( \int_0^{|\nabla u_n|}sA(x,s)ds\right)dx=O(1)\quad\mbox{as $n\ri\infty$}.$$
Using now (A2) we deduce that 
$$\intom |\nabla u_n|^{p(x)}=O(1)\quad\mbox{as $n\ri\infty$},$$
which proves our claim \eqref{roco}.

Next, we show that if $(u_n)$ is bounded in $\w10$ and satisfies for all $v\in\w10$
$$\ee'(u_n)(v)=\intom A(x,|\nabla u_n|)\nabla u_n\cdot\nabla vdx-\intom f(u_n)vdx=o(1)\quad\mbox{as $n\ri\infty$}$$
then $(u_n)$ is relatively compact. Using (f1), it is enough to show that a subsequence of $(|u_n|^{q(x)-1})$ is convergent in $W^{-1,p'(x)}(\Omega)$. By Sobolev embeddings for variable exponents, this property follows if we show that $(|u_n|^{q(x)-1})$ is relatively compact in the variable exponent Lebesgue space $L^{Np(x)/[(N+1)p(x)-N]}(\Omega)$, which is the dual space of $L^{p^*(x)}(\Omega)$.

We first observe that, up to a subsequence,
$$u_n\ri u\in L^{p^*(x)}(\Omega) \quad\mbox{a.e. as $n\ri\infty$}.$$
Fix $\delta>0$. Applying the Egorov theorem, there exists an open set $\omega\subset\Omega$ with $|\omega|<\delta$ such that
$$u_n\ri u\quad\mbox{uniformly in $\Omega\setminus\omega$.}$$
This shows that it is enough to prove that
$$\int_\omega ||u_n|^{q(x)-1}-|u|^{q(x)-1}|^{Np(x)/[(N+1)p(x)-N]}dx$$
can be made arbitrarily small.

Applying the Young inequality we find that there exists $C>0$ such that
$$\int_\omega \left(|u|^{q(x)-1} \right)^{Np(x)/[(N+1)p(x)-N]}dx\leq C\int_\omega\left( |u|^{p^*(x)}+1\right)dx,$$
which can be made as small as we wish, by choosing $\delta>0$ small enough.

Fix $\ep>0$. Since $q(x)<p^*(x)$, there exists $C_ep>0$ such that for all $n$
$$\int_\omega \left(|u_n|^{q(x)-1} \right)^{Np(x)/[(N+1)p(x)-N]}dx\leq\epsilon\int_\omega |u_n|^{p^*(x)}dx+C_\ep\,|\omega|\,.$$
Using now \eqref{roco} in combination with Sobolev embeddings for variable exponents we obtain
$$\int_\omega \left(|u_n|^{q(x)-1} \right)^{Np(x)/[(N+1)p(x)-N]}dx\leq C\epsilon+C_\ep\,|\omega|\,,$$
which can be made small enough by choosing $\delta>0$ sufficiently small.

This concludes the proof of the Palais-Smale property and of Theorem \ref{t1}. \qed

\section{Proof of Theorem \ref{t2}}

The spaces $\w10$  and $W^{-1,p'(x)}(\Omega)$ are reflexive and separable Banach spaces. Thus, by \cite{fanz}, there exist
 $\{e_j\}\subset X$ and $\{e_j^*\}\subset X^*$ such that
\begin{eqnarray*}
&&\hspace{-2cm}\w10 =\overline{\textrm{ span }\{e_j~:j=1,2,...\}},\hspace{0.5cm}W^{-1,p'(x)}(\Omega)=\overline{\textrm{ span }\{e_j^*~:j=1,2,...\}},
\end{eqnarray*}
and
$$ \langle e_i,e_j^*\rangle=\left\{
\begin{array}{ll}
1        & \mbox { if }\hspace{0.3cm} i=j,\\
0       & \mbox { if }\hspace{0.3cm}i\neq j,
\end{array}
\right.
$$
where $\langle\cdot\,,\,\cdot\rangle$ denotes the duality product between $X$ and $X^*$. We define
\begin{eqnarray*}
&&X_j=\textrm{ span}\,\{e_j\},\hspace{0.5cm}Y_k=\bigoplus_{j=1}^{k}X_j,\hspace{0.5cm}Z_k=\overline{\bigoplus_{j=k}^{\infty}X_j}.
\end{eqnarray*}

The proof of Theorem \ref{t2} is based on the following basic critical point theorem.

\begin{theorem}(Fountain Theorem, see \cite{Will}).\label{lem3.9}
Let $\ee\in C^{1}(X)$ be an even functional, where $(X,||\,.\,||)$ is a separable and reflexive Banach space. Suppose that for every $k\in \mathbb{N}$ large enogh, there exist $\rho_k>r_k>0$ such that
\begin{enumerate}
\item[$(\mathbf i)$] $\inf \{\ee(u)~:~u\in Z_k,~||u||=r_k\}\rightarrow+\infty$ as $k\rightarrow+\infty$.
\item[$(\mathbf ii)$] $\max  \{\ee(u)~:~u\in Y_k,~||u||=\rho_k\}\leq0$.
\item[$(\mathbf iii)$] $\ee$ satisfies the Palais-Smale condition for every $c>0.$
\end{enumerate}
Then $\ee$ has a sequence of critical values tending to $+\infty$.
\end{theorem}

The energy functional $\ee$ is even and satisfies the Palais-Smale condition. We show that hypotheses (i) and (ii) in the statement of Theorem \ref{lem3.9} are fulfilled. 

{\it Verification of (i)}. Fix $u\in Z_k$ with $\|u\|=r_k$, where $r_k>0$ will be specified later. Using hypotheses (A3), (A4) and (f1) we obtain
$$\begin{array}{ll}
\di\ee (u)&\di = \intom\left(\int_0^{|\nabla u(x)|}sA(x,s)ds-F(x,u)\right)dx\\
&\di \geq c_1\,|\nabla u|^{p_-}_{p(\cdot)}-c_2\intom (1+|u|^{q(x)})dx\\
&\di \geq C_1\,\| u\|^{p_-}-C_2\max\{|u|^{q_+}_{q(\cdot)},|u|^{q_-}_{q(\cdot)}\}-C_3.\end{array}$$
Assuming that $\max\{|u|^{q_+}_{q(\cdot)},|u|^{q_-}_{q(\cdot)}\}=|u|^{q_+}_{q(\cdot)}$, we obtain the estimate
$$\ee (u)\geq C_1\,\| u\|^{p_-}-C_2\,|u|^{q_+}_{q(\cdot)}-C_3.$$

 Denote
$$\alpha_{k}=\sup\{|u|_{L^{q(x)}(\Omega)};~||u||=1,~u\in Z_k\}.$$
Then  $\lim_{k\rightarrow\infty}\alpha_k=0$, see \cite{Fan6}.

Since $u\in Z_k$, we deduce that
$$\ee (u)\geq C_1\,\| u\|^{p_-}-C_2\,\alpha_k^{q_+}\|u\|^{q_+}-C_3.$$
Since $p_-<q_+$ and $\alpha_k\ri 0$, it follows that
$$r_k:=\left(\frac{C_2}{C_1}\alpha_k^{p_+}\right)^{1/(p_--q_+)}\ri +\infty\quad\mbox{as $k\ri\infty$}.$$
Taking $\|u\|=r_k$ we conclude that $\ee (u)\ri+\infty$, hence (i) is fulfilled.

\smallskip
{\it Verification of (ii)}. Fix $u\in Y_k$ with $\|u\|=1$ and let $\rho_k$ be a positive number, which will be defined later. We have seen in relation \eqref{geo4} that there exists $\alpha>p_+$ such that
$$\ee (\rho_ku)\leq C_1(u)\,\rho_k^{p_+}-C_2\,\rho_k^\alpha+C_3,$$
where $C_1$, $C_2$, $C_3$ are positive constants. Taking $\rho_k>r_k$ and using the fact that $\alpha>p_+$, we deduce that $\ee(\rho_ku)\ri -\infty$ as $k\ri\infty$. This implies that
$$\max\{\ee (u);\ u\in Y_k,\ \|u\|=\rho_k\}\leq 0,$$
for every $\rho_k$ large enough. Applying the Fountain Theorem, we complete the proof of Theorem \ref{t2}. \qed

\section{Proof of Theorem \ref{t3}}
For all $u\in\w10$ we define
$$I(u):=\intom\left(\int_0^{|\nabla u(x)}sA(x,s)ds\right)dx\quad\mbox{and}\quad J(u):=\intom F(x,u)dx.$$

It follows that
$$\lambda_1=\inf_{u\in\w10\setminus\{0\}}\frac{I(u)}{J(u)}\,.$$

Let $(u_n)\subset\w10$ such that $u_n\rightharpoonup u$. We prove that
\begin{equation}\label{auu}I(u)\leq\liminf_{n\ri\infty}I(u_n)\quad\mbox{and}\quad J(u)=\lim_{n\ri\infty}J(u_n).\end{equation}

Since $\ee_0'$ is monotone, it follows that $I$ is convex. Therefore
$$I(u_n)\geq I(u)+I'(u)(u_n-u),$$
hence $\liminf_{n\ri\infty}I(u_n)\geq I(u)$.

Next, using assumptions (f1) and (f3) we deduce that there is a positive constant $C$ such that for all $(x,u)\in\Omega\times\RR$
\begin{equation}\label{lju}|F(x,u)|\leq |u|^{p_+}+C\, |u|^{q(x)}.\end{equation}
Since $p(x),\,q(x)<p^*(x)$ and $(u_n)$ is bounded in $\w10$ we can assume that, up to a subsequence, $u_n\ri u$ both in $L^{p_+}(\Omega)$ and in $L^{q(x)}(\Omega)$. Using relation \eqref{lju} we deduce that $J(u_n)\ri J(u)$ as $n\ri\infty$.

Our hypotheses imply that for all $u\in\w10$ with $\|u\|$ small enough we have
$$I(u)\geq C_1\,\|u\|^{p_+}$$
and
$$\begin{array}{ll}\di 0&\di\leq J(u)=\intom F(x,u)dx\leq C\intom\frac{|u|^{q(x)}}{q(x)}dx\\
&\di \leq C_2|u|^{q_+}_{q(\cdot)}+C_3|u|^{q_-}_{q(\cdot)}\leq C_2\|u\|^{q_+}+C_3\|u\|^{q_-}\,.\end{array}$$
Since $p_+<q_-$ we deduce that
\begin{equation}\label{cozero}
\lim_{\|u\|\ri 0}\frac{I(u)}{J(u)}=+\infty\,.
\end{equation}
Using hypothesis (A5) and similar energy estimates as above we deduce that
\begin{equation}\label{coinfty}
\lim_{\|u\|\ri \infty}\frac{I(u)}{J(u)}=+\infty\,.
\end{equation}

\smallskip
{\it Step 1.} We prove that $\lambda_1>0$.

Our assumptions imply that $\lambda_1\geq 0$. Arguing by contradiction and supposing that $\lambda_1=0$, we find a sequence $(u_n)\subset\w10\setminus\{0\}$ such that
\begin{equation}\label{oxo2}
\lim_{n\ri \infty}\frac{I(u_n)}{J(u_n)}=0\,.
\end{equation}
We have already remarked that
\begin{equation}\label{oxo3}
\frac{I(u)}{J(u)}\geq\frac{C_1\,\min\{\|u\|^{p_+},\|u\|^{p_-}\}}{C_2\,\|u\|^{q_+}+C_3\,\|u\|^{q_-}}\,.
\end{equation}
Since $p_+<q_-$, relations \eqref{oxo2} and \eqref{oxo3} imply that $(u_n)$ is unbounded. Using now \eqref{coinfty} we contradict our assumption \eqref{oxo2}.

\smallskip
{\it Step 2.} We show that problem \eqref{problem1} has a solution for $\lambda=\lambda_1$.

Let $(u_n)\subset\w10\setminus\{0\}$ be such that 
\begin{equation}\label{ola}\lim_{n\ri \infty}\frac{I(u_n)}{J(u_n)}=]\lambda_1\,.\end{equation}
Using \eqref{auu} we deduce that $(u_n)$ is bounded. Thus, up to a subsequence, 
$$u_n\rightharpoonup u\quad\mbox{in $\w10$}.$$
Assuming now that $u=0$, relation \eqref{auu} shows that $J(u_n)\ri 0$. Thus, by \eqref{ola}, we also have $I(u_n)\ri 0$. But for all $u\in\w10$
$$I(u)\geq C_1\,\min\{\|u\|^{p_+},\|u\|^{p_-}\}.$$
We deduce that $\|u_n\|\ri 0$. Using now relation \eqref{coinfty} we obtain a contradiction. This proves that $u\not=0$. Using now \eqref{auu} we conclude that $I(u)=\lambda_1 J(u)$, hence $\lambda_1$ is an eigenvalue of problem \eqref{problem1}.

\smallskip
{\it Step 3.} Every $\lambda>\lambda_1$ is an eigenvalue of problem \eqref{problem1}.

Fix $\lambda>\lambda_1$. The energy functional associated to problem \eqref{problem1} is 
$${\mathcal E}_\lambda (u):=I(u)-\lambda J(u).$$
Relation \eqref{coinfty} shows that $\ee_\lambda$ is coercive, that is, $\ee_\lambda(u)\ri +\infty$ as $\|u\|\ri\infty$. But, by \eqref{auu}, $\ee_\lambda$ is lower semi-continuous, hence it has a global minimizer $w\in\w10$. On the other hand, since $\lambda>\lambda_1$, there exists $v\in\w10$ such that
$$\lambda_1<\frac{I(v)}{J(v)}<\lambda\,.$$
This shows that $\ee_\lambda (v)<0$, so $\ee_\lambda (w)<0$. We conclude that $w\not=0$ and $w$ is a critical point of $\ee_\lambda$, hence a nontrivial solution of problem \eqref{problem1}. \qed 

\medskip
 \indent {\bf Acknowledgements.}  
 The author acknowledges the support by the Slovenian Research Agency grants
 P1-0292, J1-4144, and J1-5435.

\end{document}